\documentclass[11pt,a4paper]{article}

\usepackage{amsmath,amssymb,amsthm}

\theoremstyle{plain}
\newtheorem{satz}{Theorem}[section]
\newtheorem*{namedtheorem}{\theoremname}
\newcommand{\theoremname}{testing}
\newenvironment{satzmitname}[1]{\renewcommand{\theoremname}{#1}
   \begin{namedtheorem}}
   {\end{namedtheorem}}
\newtheorem{lemma}[satz]{Lemma}
\newtheorem{prop}[satz]{Proposition}
\newtheorem{cor}[satz]{Corollary}

\newenvironment{bew}{\begin{proof}}{\end{proof}}
\theoremstyle{definition}
\newtheorem*{defn}{Definition}

\newcommand{\Q}{\ensuremath{\mathbb Q}}
\newcommand{\Z}{\ensuremath{\mathbb Z}}
\newcommand{\N}{\ensuremath{\mathbb N}}
\newcommand{\R}{\ensuremath{\mathbb R}}
\newcommand{\C}{\ensuremath{\mathbb C}}

\newcommand{\tp}{^\mathrm T}
\newcommand{\tT}{\mathrm T}
\newcommand{\tr}{\ensuremath{\operatorname{tr}}}
\newcommand{\Gal}{\ensuremath{\operatorname{Gal}}}

\newcommand{\M}{\ensuremath{\mathrm{M}}}
\newcommand{\GL}{\ensuremath{\mathrm{GL}}}
\newcommand{\SL}{\ensuremath{\mathrm{SL}}}
\newcommand{\SP}{\ensuremath{\mathrm{Sp}}}
\newcommand{\U}{\ensuremath{\mathrm U}}

\newcommand{\quer}{\overline}
\newcommand{\wtilde}{\widetilde}
\newcommand{\de}{\mathfrak}

\renewcommand{\leq}{\leqslant}
\renewcommand{\epsilon}{\varepsilon}

\begin{document}

\title{Conjugacy classes of \hbox{$p$}-torsion in symplectic
        groups over \hbox{$S$}-integers}
\author{Cornelia Minette Busch\thanks{Research supported by a grant
"Estancias de j\'ovenes doctores y tecn\'ologos extranjeros en Espa\~na"
(SB 2001 - 0138) from the Ministerio de Educaci\'on,
Cultura y Deporte.}
}
\date{September 2005}
\maketitle

\begin{abstract}
For any odd prime $p$ we consider representations of a group of order
$p$ in the symplectic group $\SP(p-1,\Z[1/n])$ of
$(p-1)\/\times\/(p-1)$-matrices over the ring $\Z[1/n]$, $0\neq n\in\N$.
We construct a relation between the conjugacy classes of subgroups $P$ of
order $p$ in the symplectic group and the ideal class group in the ring
$\Z[1/n]$.
This is used for the study of these classes. In particular we determine
the centralizer $C(P)$ and $N(P)/C(P)$ where $N(P)$ denotes the normalizer.
\end{abstract}

{\small{\hspace*{2mm}

\hspace*{2mm}
}}


\section{Introduction}

We define the group of symplectic matrices $\SP(2n,R)$ over a
ring $R$ to be the subgroup of matrices $M\in\GL(2n,R)$ that satisfy
$$M\tp JM =J=\begin{pmatrix}\phantom{-}0 & 1 \\ -1 & 0\end{pmatrix}$$
where $1\in\M(n,R)$ denotes the identity.
Our motivation for studying subgroups of odd prime order $p$ in
the symplectic group $\SP(p-1,\Z[1/n])$, $0\neq n\in\N$, is given by the
fact that the $p$-primary part of the Farrell cohomology of
$\SP(p-1,\Z[1/n])$ is determined by the Farrell cohomology of the
normalizer of subgroups of order $p$ in $\SP(p-1,\Z[1/n])$
(see Brown~\cite{brownb}). First we consider
the conjugacy classes of elements of order $p$ in $\SP(p-1,\Z[1/n])$ and
get the following result.

\pagebreak

\begin{satzmitname}{Theorem~\ref{Konj'kl von Elem. in Sp(p-1,Z[1/n])}}
The number of conjugacy classes of matrices of order $p$ in
$\SP(p-1,\Z[1/n])$, $0\neq n\in\Z$, is
$$|\mathcal C_0| 2^{\frac{p-1}{2}+\tau}$$
where $\mathcal C_0$ is the ideal class group of $\Z[1/n][\xi]$ and
$\tau$ is the number of inert primes in $\Z[\xi +\xi^{-1}]$ that lie over
primes in $\Z$ that divide $n$.
\end{satzmitname}

In order to prove this theorem we establish a relation between
some ideal classes in $\Z[1/n][\xi]$, $\xi$ a primitive $p$th root
of unity, and the conjugacy classes of matrices of order $p$.
We define equivalence classes
$[\de a,a]$ of pairs $(\de a,a)$ where $\de a\subseteq \Z[1/n][\xi]$ is
an ideal with $\de a\quer{\de a}=(a)$ and the equivalence relation is
$$
\begin{array}{rcl}
(\de a,a)\sim (\de b,b)
&
\Leftrightarrow
&
\exists \lambda,\mu\in\Z[1/n][\xi]\setminus \{0\}
\\
&
&
\lambda\de a =\mu\de b,\ \lambda\quer\lambda a=\mu\quer\mu  b\,.
\end{array}
$$
We show that a bijection exists between the conjugacy classes of elements
of order $p$ in $\SP(p-1,\Z[1/n])$ and the set of equivalence classes
$[\de a,a]$. Sjerve and Yang (see~\cite{sjerve}) construct an
analogous bijection for $\SP(p-1,\Z)$. We use the bijection
described above in order to study the subgroups of order $p$ in
$\SP(p-1,\Z[1/n])$. We consider the case where $n\in\Z$ is such that
$\Z[1/n][\xi]$ and $\Z[1/n][\xi +\xi^{-1}]$ are principal ideal domains
because in this case the ideal class group of those rings is trivial.
We get the following results.

\begin{satzmitname}{Theorem~\ref{C in Sp(p-1,Z[1/n])}}
Let $n\in\Z$ be such that $\Z[1/n][\xi]$ and $\Z[1/n][\xi +\xi^{-1}]$ are
principal ideal domains. Then the centralizer $C(P)$ of a
subgroup $P$ of order $p$ in $\SP(p-1,\Z[1/n])$ is
$$C(P)\cong \Z/2p\Z\times\Z^{\sigma^+}$$
where $\sigma^+ =\sigma$ if $p\nmid n$, $\sigma^+=\sigma +1$ if
$p\mid n$ and $\sigma$ is the number of primes in $\Z[\xi +\xi^{-1}]$ that
split in $\Z[\xi]$ and lie over primes in $\Z$ that divide $n$.
\end{satzmitname}

\begin{satzmitname}{Theorem~\ref{N/C in Sp(p-1,Z[1/n])}}
Let $n\in\Z$ be such that $\Z[1/n][\xi]$ and $\Z[1/n][\xi +\xi^{-1}]$
are principal ideal domains and moreover $p\mid n$. Let $N(P)$ denote
the normalizer and $C(P)$ the centralizer of a subgroup $P$ of order
$p$ in $\SP(p-1,\Z[1/n])$. Then
$$N(P)/C(P)\cong\Z/j\Z$$
where $j\mid p-1$, $j$ odd. For each $j$ with $j\mid p-1$, $j$ odd,
exists a subgroup of order $p$ in $\SP(p-1,\Z[1/n])$
with $N(P)/C(P)\cong\Z/j\Z$.
\end{satzmitname}

An application of these theorems is given in \cite{buschperFar};
moreover they are a generalization of the results of
Naffah~\cite{nadim} on the normalizer of $\SL(2,\Z[1/n])$.

Let
$\U\bigl((p-1)/2\bigr)\subset\GL\bigl((p-1)/2,\C\bigr)$
be the group of unitary matrices. We consider the homomorphism
$$
\begin{array}{rcl}
\U\bigl(\frac{p-1}{2}\bigr) & \longrightarrow & \SP(p-1,\R) \\
X = A+ iB & \longmapsto & \begin{pmatrix}
                        \phantom{-}A & B \\ -B & A
                      \end{pmatrix}\,
\end{array}
$$
where $A,B\in\M(n,\R)$. In \cite{BuschSC} a condition is given
for the matrix $X$ such that the image of $X$ is conjugate to
a matrix of order $p$ in $\SP(p-1,\Z)$. This is used in
\cite{BuschHFSp} to analyze the subgroups of order $p$ in
$\SP(p-1,\Z)$ by considering the corresponding subgroups in
$\U((p-1)/2)$. Here we avoid the unitary group by taking an
arithmetical approach.


\section{A recall of algebraic number theory}


For the convenience of the reader, we give a short introduction to
algebraic number theory. More details and the proofs
can be found in the books of Lang~\cite{langant}, Neukirch~\cite{neukirch}
and Washington~\cite{wash}.

Let $p$ be an odd prime and let $\xi$ be a primitive $p$th root
of unity. Then $\Z[\xi]$ is the ring of integers of the cyclotomic field
$\Q(\xi)$ and $\Z[\xi +\xi^{-1}]$ is the ring of integers of the maximal real
subfield $\Q(\xi +\xi^{-1})$ of $\Q(\xi)$.
For an integer $0\neq n\in\Z$ we consider the ring $\Z[1/n]$ and
the extensions $\Z[1/n][\xi]$ and $\Z[1/n][\xi +\xi^{-1}]$. It is well-known that
$\Z[1/n][\xi]$ and $\Z[1/n][\xi +\xi^{-1}]$ are Dedekind rings.
For $j=1,\dots,p-1$ let the Galois automorphism
$\gamma_j\in\Gal(\Q(\xi)/\Q)$ be given by $\gamma_j(\xi)=\xi^j$.
To simplify the notations, we define $x^{(j)}:=\gamma_j(x)$ for any
$x\in\Q(\xi)$ and $\gamma_j\in\Gal(\Q(\xi)/\Q)$ as above.
The Galois automorphism $\gamma_j$ acts componentwise on a vector
in $\Q(\xi)^k$.

Let $A$ be a Dedekind ring and $K$ the quotient field of $A$. Let $L$ be
a finite separable extension of $K$ and $B$ the integral closure
of $A$ in $L$. Let $\de a$ be an additive subgroup of $L$.
The complementary set $\de a'$ of $\de a$ is
the set of $x\in L$ such that $\tr_{L/K}(x\de a)\subseteq A$.
The different of the extension $B/A$ is defined to be
$$D_{B/A}:={B'}_{L/K}^{-1}.$$
In $\Z[\xi]$ the different is generated by $D=p\xi^{(p+1)/2}/(\xi -1)$.
It is a principal ideal. This is also true for $\Z[1/n][\xi]$
(see Lang~\cite{langant} or Serre~\cite{serrelf}).

Let $\mathcal O$ be the ring of integers of a number field $K$. Let
$G=\Gal(K/\Q)$ be the Galois group of the extension and let
$\de q$ be a prime ideal of $\mathcal O$. The subgroup
$$G_{\de q}=\{\sigma\in G\mid \sigma\de q=\de q\}$$
is called the decomposition group of $\de q$ over $\Q$. The fixed field
$$Z_{\de q}=\{x\in K\mid \sigma x=x \text{ for all }\sigma\in G_{\de q}\}$$
is called the decomposition field of $\de q$ over $\Q$.
The decomposition group of a prime ideal $\sigma\de q$ that is conjugate
to $\de q$ is the conjugate subgroup
$G_{\sigma\de q} = \sigma G_{\de q}\sigma^{-1}$.
Let $\de q\subset \mathcal O$ be a prime ideal in $\mathcal O$ over the prime
$(q)$ in $\Z$. Let $\kappa(\de q):=\mathcal O/\de q$ and $\kappa(q):=\Z/q\Z$.
The degree $f_{\de q}$ of the extension of fields
$\kappa(\de q)/\kappa(q)$ is called the residue class degree
of $\de q$. We recall the following property.
For any prime $q\neq p$ let $f_q\in\N$ be the smallest positive integer
such that
$$q^{f_q} \equiv 1\mod p.$$
Then $(q)=(\de q_1\cdots\de q_r)$
where $\de q_1,\dots,\de q_r$ are pairwise different prime ideals
in $\Q(\xi)$ and all have residue class degree $f_q$ (see
Neukirch~\cite{neukirch}).

Let $p$, $q$ and $\xi$ be as above.
Let $\de q^+\subseteq \Z[\xi +\xi^{-1}]$ be a prime ideal that lies over $q$.
We consider the ideal $\de q^+\Z[\xi]\subset\Z[\xi]$ generated by $\de q^+$.
Any prime $q\neq p$ is unramified and the prime $p$ ramifies.
Let $\sigma\in G:=\Gal(\Q(\xi)/\Q)$ with $\sigma (x)=\quer x$.
The Galois group $G$ acts transitively on the set of prime ideals over
$q$. It is known that $f_q=|G_{\de q}|$. We have the following three cases.
\newline
The prime $\de q^+$ is inert: $\de q^+\Z[\xi]=\de q$, a prime ideal in
$\Z[\xi]$ that lies over $q$.
\begin{center}
 \begin{minipage}[t]{11.5cm}
  \begin{tabular}{lll}
  \makebox[1.3cm][r]{$\de q^+\Z[\xi]=\de q$}
    & $\Leftrightarrow$ &
     \begin{minipage}[t]{8cm}$ \de q=\quer{\de q}$
     \end{minipage}\\
    & $\Leftrightarrow$ &
     \begin{minipage}[t]{8cm} $\sigma\in G_{\de q}$, i.e. $G_{\de q}$
     contains an element of order $2$,
     \end{minipage}\\
    & $\Leftrightarrow$ &
     \begin{minipage}[t]{8cm} $f_q$ is even.
     \end{minipage}
  \end{tabular}
 \end{minipage}
\end{center}
Primes that split in $\Z[\xi]$: $\de q^+\Z[\xi]=\de q\quer{\de q}$
where $\de q$ is a prime ideal in $\Z[\xi]$ that lies over $q$.
\begin{center}
 \begin{minipage}[t]{11.5cm}
  \begin{tabular}{lll}
  \makebox[1.3cm][r]{$\de q^+\Z[\xi]=\de q\quer{\de q}$}
    & $\Leftrightarrow$ &
     \begin{minipage}[t]{8cm} $\de q\neq\quer{\de q}$ \end{minipage}\\
    & $\Leftrightarrow$ &
     \begin{minipage}[t]{8cm} $\sigma\not\in G_{\de q}$, i.e.
     $G_{\de q}$ does not contain an element of order $2$,
     \end{minipage}\\
    & $\Leftrightarrow$ &
     \begin{minipage}[t]{8cm} $f_q$ is odd.
     \end{minipage}
  \end{tabular}
 \end{minipage}
\end{center}
The ramified case:  $\de p^+\Z[\xi]=\de p^2$ where $\de p:=(1-\xi)$ is
the only prime ideal in $\Z[\xi]$ that lies over $p$. Moreover
$\de p^+\Z[\xi] :=((1-\xi)(1-\xi^{-1}))=\de p \quer{\de p}$
is the only prime ideal in $\Z[\xi +\xi^{-1}]$ that lies over $p$.

Let $\mathcal O_K$ be a Dedekind ring and let $S$ be a finite set of prime
ideals $\de q\subseteq\mathcal O_K$. We define
$$\mathcal O_K^S:=
  \left\{\left.
   \frac{f}{g}\ \right|\
    f,g\in\mathcal O,\ g\not\equiv 0\mod\de q\text{ for }\de q\not\in S
  \right\}.$$
Let $K$ be the quotient field of $\mathcal O_K$.
We call the group $(\mathcal O_K^S)^*$ the group of $S$-units of $K$.
Let $\mathcal C(\mathcal O_K)$, resp.\ $\mathcal C(\mathcal O_K^S)$,
denote the ideal class group of $\mathcal O_K$, resp.\ $\mathcal O_K^S$.

\begin{prop}\label{Struktur S-units}
For the group $(\mathcal O_K^S)^*$ defined above we have an isomorphism
$$(\mathcal O_K^S)^*\cong \mu(K)\times \Z^{|S|+r+s-1}$$
where $\mu(K)$ denotes the group of roots of unity of $K$, $r$ denotes
the number of real embeddings of $K$ and $s$ denotes the number of
conjugate pairs of complex embeddings of $K$.
\end{prop}

\begin{bew}
See Neukirch~\cite{neukirch}.
\end{bew}

Therefore
\begin{align*}
(\mathcal O_K^S)^* &\cong \mu(K)\times \Z^{r+s-1}\times \Z^{|S|} \\
   &\cong \mathcal O_K^* \times \Z^{|S|}.
\end{align*}
%


\section{Matrices of order $p$}


\subsection{A relation between matrices and ideal classes}

The results obtained in this section are based on the bijection given
by Proposition~\ref{matrizen ideale}. Sjerve and Yang prove in
\cite{sjerve} the analogous statement of this proposition for the group
$\SP(p-1,\Z)$. Since for our purpose it is important to understand
the bijection and some proofs need a slightly different approach
for the group $\SP(p-1,\Z[1/n])$, we present in this subsection
some of the proofs for the convenience of the reader.

\begin{defn}
Let $I$ be the set of pairs $(\de a,a)$\, where
$\de a\subseteq\Z[1/n][\xi]$ is a $\Z[1/n][\xi]$-ideal and
$0\neq a\in\Z[1/n][\xi]$ is such that
$\de a\overline{\de a}=(a)\subseteq\Z[1/n][\xi]$. Here
$\overline{\de a}$ denotes the ideal generated by the complex
conjugate of the elements of $\de a$.
We define an equivalence relation on $I$.
$$\begin{array}{rcl}
(\de a,a)\sim (\de b,b) & \Leftrightarrow &
    \exists \lambda,\mu\in\Z[1/n][\xi],\ \lambda,\mu\neq 0 \\
& & \lambda\de a = \mu\de b,
    \ \lambda\overline{\lambda}a=\mu\overline{\mu}b. \\
\end{array}$$
Let $[\de a,a]$ denote the equivalence class of the pair
$(\de a,a)$ and let $\mathcal I$ be the set of equivalence
classes $[\de a,a]$.
\end{defn}

\begin{lemma}\label{B pairs equal P}
Let $(\de a,a)$ be a pair consisting of a $\Z[1/n][\xi]$-ideal
$\de a\subseteq\Z[1/n][\xi]$ and $0\neq a\in\Z[1/n][\xi]$.
Then $(\de a,a)\in I$ if and only if a $\Z[1/n]$-basis
$\alpha_1,\dots,\alpha_{p-1}$ of $\de a$ exists such that
$$\alpha\tp J \quer\alpha^{(i)}=\delta_{1i}aD$$
where $D=p\xi^{(p+1)/2}/(\xi -1)$ and
$\alpha=(\alpha_1,\dots,\alpha_{p-1})\tp$.
\end{lemma}

\begin{bew}
The proof is analogous to the proof of Lemma 2.3 in \cite{sjerve}.
\end{bew}

\begin{lemma}\label{matrix gleich}
Let $M$, $N$ be two $(p-1)$\/$\times$\/$(p-1)$-matrices over $\Z[1/n]$ and
let
$$\alpha=(\alpha_1,\dots,\alpha_{p-1})\tp\in\Z[1/n][\xi]^{p-1}$$
where $\alpha_1,\dots,\alpha_{p-1}$ are $\Z[1/n]$-linear independent.
If for $i=1,\dots,p-1$
$$\alpha\tp M \quer\alpha^{(i)}=\alpha\tp N \quer\alpha^{(i)}$$
then we have $M=N$.
\end{lemma}

\begin{bew}
It suffices to prove the case $N=0$ because
$$\alpha\tp M \quer\alpha^{(i)}
    = \alpha\tp N \quer\alpha^{(i)}
  \quad\Leftrightarrow \quad
  \alpha\tp (M-N) \quer\alpha^{(i)}
    = 0 = \alpha\tp 0 \quer\alpha^{(i)}.$$
Let $a_i=\alpha\tp M \quer\alpha^{(i)}$, then
$a_i^{(k)}=\alpha^{(k)\tT}M (\quer\alpha^{(i)})^{(k)}$.
For all $k,l$ with $1\leq k,l\leq p-1$ let $i$ be such that
$1\leq i\leq p-1$ and $ki\equiv l\mod p$. Then
$(\quer\alpha^{(i)})^{(k)}=\quer\alpha^{(l)}$ and therefore
$\alpha^{(k)\tT} M \quer\alpha^{(l)} = 0$ for $k,l=1,\dots,p-1$.
This implies $A\tp MB=0$ where
$$A:=\bigl(\alpha_i^{(j)}\bigr) \text{ and }
  B:=\bigl(\quer\alpha_i^{(j)}\bigr)$$
are $(p-1)$\/$\times$\/$(p-1)$-matrices.
Since $\alpha_1,\dots,\alpha_{p-1}$ are $\Z[1/n]$-linear independent
we have $\det A\neq 0$ and $\det B\neq 0$. But this yields $M=0$.
\end{bew}

\begin{prop}\label{matrizen ideale}
A bijection $\psi$ exists between the set of conjugacy classes
of elements of order $p$ in $\SP(p-1,\Z[1/n])$ and the set of
equivalence classes of pairs $[\de a,a]\in \mathcal I$.
\end{prop}

In order to prove this proposition, we first construct the
bijection and then we show that the mapping we constructed is a
bijection (Lemma~\ref{psi injektiv}, Lemma~\ref{psi surjektiv}).

Let $Y\in\SP(p-1,\Z[1/n])$ be of order $p$. The eigenvalues of $Y$
are the primitive $p$th roots of unity. An eigenvector
$$\alpha =(\alpha_1,\dots,\alpha_{p-1})\tp \in (\Z[1/n])^{p-1}$$
exists for the eigenvalue $\xi=e^{i2\pi/p}$, i.e. $Y\alpha =\xi\alpha$.
The $\alpha_1,\dots,\alpha_{p-1}$ are $\Z[1/n]$-linear independent.
Let $\de a$ be the $\Z[1/n]$-module generated by
$\alpha_1,\dots,\alpha_{p-1}$.
Let $a=D^{-1}\alpha\tp J\quer\alpha$.
Then $\de a\subseteq\Z[1/n][\xi]$ is a
$\Z[1/n][\xi]$-ideal and $a=\quer a$.

\begin{lemma}
The pair $(\de a,a)$ we construct above is an element of $I$.
\end{lemma}

\begin{bew}
Because of Lemma~\ref{B pairs equal P}
it suffices to show that $\alpha\tp J\quer\alpha^{(i)}=0$ for
$i=2,\dots,p-1$. Since $Y\alpha=\xi\alpha$ we have
$$Y\alpha^{(i)} = \xi^i \alpha^{(i)} \text{ and }
  Y\quer\alpha^{(i)} = \frac{1}{\xi^i}\, \quer\alpha^{(i)},$$
$2\leq i\leq p-1$. Therefore
$$\alpha\tp J\quer\alpha^{(i)} =
  \frac{\xi^i}{\xi}\,\alpha\tp Y\tp JY\quer\alpha^{(i)} =
  \frac{\xi^i}{\xi}\,\alpha\tp J \quer\alpha^{(i)}$$
where the last equation follows from the fact that
$Y\in\SP(p-1,\Z[1/n])$. Since $\xi\neq \xi^j$ we get
$\alpha\tp J\quer\alpha^{(i)}=0$.
\end{bew}

Let $Y,\wtilde Y\in\SP(p-1,\Z[1/n])$ be matrices of odd prime order $p$.
Let $\alpha\in (\Z[1/n][\xi])^{p-1}$, resp.\ $\beta\in (\Z[1/n][\xi])^{p-1}$
be an eigenvector of $Y$, resp.\ $\wtilde Y$, to the eigenvalue $\xi$, i.e.
$Y\alpha=\xi\alpha$ and $\wtilde Y\beta=\xi\beta$. Let
$\alpha=(\alpha_1,\dots,\alpha_{p-1})\tp$,
$\beta=(\beta_1,\dots,\beta_{p-1})\tp$.
Let $\de a\subseteq\Z[1/n][\xi]$, resp.\ $\de b\subseteq\Z[1/n][\xi]$,
be the ideal with $\Z[1/n]$-basis $\alpha_1,\dots,\alpha_{p-1}$,
resp.\ $\beta_1,\dots,\beta_{p-1}$. We define
$a=D^{-1}\alpha\tp J\quer\alpha$ and $b=D^{-1}\beta\tp J\quer\beta$.
We show the injectivity of $\psi$.

\begin{lemma}\label{psi injektiv}
Let $Y,\wtilde Y\in\SP(p-1,\Z[1/n])$ be matrices of odd prime order $p$.
Then $Y$ and $\wtilde Y$ are conjugate if and only if
$[\de a,a]=[\de b,b]$.
\end{lemma}

\begin{bew}
Let $Y$ and $\wtilde Y$ be conjugate. Then $Q\in\SP(p-1,\Z[1/n])$ exists
such that $\wtilde Y = Q^{-1}YQ$. Then $Q\wtilde Y=YQ$ and for the
eigenvector $\beta$ to the eigenvalue $\xi$ of $\wtilde Y$ we get
$$YQ\beta = Q\wtilde Y\beta = \xi Q\beta$$
and therefore $Q\beta$ is an eigenvector of $Y$. But $\alpha$ is also an
eigenvector to the eigenvalue $\xi$ of $Y$. So $\lambda,\mu\in\Z[1/n][\xi]$,
$\lambda,\mu\neq 0$, exist such that
$$\lambda\alpha =\mu Q\beta =Q\mu\beta.$$
Then $\lambda\de a=\mu\de b$ and for $a=D^{-1}\alpha\tp J\quer\alpha$,
$b=D^{-1}\beta\tp J\quer\beta$ we get
$\lambda\quer\lambda a = \mu\quer\mu \,b$.
This shows that $[\de a,a]=[\de b,b]$.

In order to show the other direction we assume that
$\lambda,\mu\in\Z[1/n][\xi]$, $\lambda,\mu\neq 0$, exist such that
$\lambda\de a = \mu\de b$ and
$\lambda\quer\lambda a =\mu\quer\mu b$. Then a matrix $Q\in\GL(p-1,\Z[1/n])$
exists such that $\lambda\alpha =\mu Q\beta$.
We have
$$
\mu Q\wtilde Y\beta
= \mu Q\xi\beta = \xi\mu Q\beta = \xi\lambda\alpha = \lambda Y\alpha
= \mu YQ\beta
$$
and therefore
$$Q\wtilde Y\beta = YQ\beta.$$
Since $\beta_1,\dots,\beta_{p-1}$ are $\Z[1/n]$-linear independent, we
have $Q\wtilde Y = YQ$ and herewith
$$\wtilde Y = Q^{-1} YQ.$$
It remains to show that $Q\in\SP(p-1,\Z[1/n])$. For $i=2,\dots,p-1$
we have
$$\beta\tp Q\tp JQ\,\quer\beta^{(i)}
  = \frac{\lambda\quer\lambda^{(i)}}{\mu\quer\mu^{(i)}} \,
     \alpha\tp J\,\quer\alpha^{(i)}
  = 0 = \beta\tp J\,\quer\beta^{(i)}$$
and for $i=1$ we have
$$\beta\tp Q\tp JQ\,\quer\beta
  = \frac{\lambda\quer\lambda}{\mu\quer\mu} \, \alpha\tp J\,\quer\alpha
  = \frac{b}{a}\,\alpha\tp J\,\quer\alpha
  = \beta\tp J\,\quer\beta$$
because $\lambda\quer\lambda\,a = \mu\quer\mu\,b$ implies that
$\frac{\lambda\quer\lambda}{\mu\quer\mu} = \frac{b}{a}$. Now it follows
from Lemma~\ref{matrix gleich} that $Q\tp JQ = J$ and this means that
$Q\in\SP(p-1,\Z[1/n])$.
\end{bew}

\begin{lemma}\label{psi surjektiv}
The mapping $\psi$ is surjective.
\end{lemma}

\begin{bew}
Let $(\de a,a)$ and $\alpha=(\alpha_1,\dots,\alpha_{p-1})\tp$ be be
like in Lemma~\ref{B pairs equal P}. Then $\xi\alpha_1,\dots,\xi\alpha_{p-1}$
is a new basis of $\de a$. Therefore $X\in\GL(p-1,\Z[1/n])$ exists
with $X\alpha=\xi\alpha$. It is evident that the order of $X$ is $p$.
We show that $X\in\SP(p-1,\Z[1/n])$. We have
$$\alpha\tp X\tp JX\quer\alpha^{(i)} =
  \frac{\xi}{\xi^i}\,\alpha\tp J\quer\alpha^{(i)} =
  \delta_{1i}\alpha\tp J\,\quer\alpha$$
hence
$$\alpha\tp X\tp JX\,\quer\alpha^{(i)} = \alpha\tp J\,\quer\alpha^{(i)}.$$
The last equation and Lemma~\ref{matrix gleich} imply that
$X\tp JX = J$ and therefore $X\in\SP(p-1,\Z[1/n])$.
\end{bew}

Let $\mathcal I$ be the set of equivalence classes of pairs
$(\de a,a)\in I$ defined above. We define a multiplication
on $\mathcal I$ by
$$[\de a,a]\cdot[\de b,b]=[\de a\de b,ab].$$
The unit is $[\Z[1/n][\xi],1]$ and the inverse of $[\de a,a]$
is $[\quer{\de a},a]$ since
$$[\de a,a]\cdot[\quer{\de a},a]=[(a),a^2]=[\Z[1/n][\xi],1].$$

\begin{lemma}
Let $(\de a,a)\in I$, $\lambda\in\Z[1/n][\xi]$, $\lambda\neq 0$. Then
\begin{enumerate}
\item $(\lambda\de a,\lambda\quer\lambda a)\in I$,
\item $(\de a,\lambda a)\in I$ if and only if
      $\lambda\in\Z[1/n][\xi +\xi^{-1}]^*$.
\end{enumerate}
\end{lemma}

\begin{bew}
Trivial.
\end{bew}

Let
$$N:\Q(\xi)\longrightarrow \Q(\xi +\xi^{-1})$$
be the norm mapping, i.e.\ $N(x)=x\quer x$ for $x\in\Q(\xi)$.
Then
$$N(\Z[1/n][\xi]^*)\subsetneq\Z[1/n][\xi +\xi^{-1}]^*.$$

\begin{lemma}\label{delta injektiv}
Let $(\de a,a),(\de a,b)\in I$. Then
$[\de a,a]=[\de a,b]$ if and only if
$$\frac{a}{b}\in N(\Z[1/n][\xi]^*).$$
\end{lemma}

\begin{bew}
Suppose that $[\de a,a]=[\de a,b]$. Then
$\lambda,\mu\in\Z[1/n][\xi]$, $\lambda,\mu\neq 0$, exist such that
$\lambda\de a=\mu\de a$ and
$\lambda\quer\lambda \,a=\mu\quer\mu \,b$.
Let $u=\mu/\lambda$, then $u\in\Z[1/n][\xi]^*$
(since $\de a = (\mu /\lambda)\,\de a$)
and $a/b=\mu\quer\mu/\lambda\quer\lambda=u\quer u$. This
shows that $a/b\in N(\Z[1/n][\xi]^*)$.
Now let $a/b=u\quer u$ for some $u\in\Z[1/n][\xi]^*$. Then
$[\de a,a]=[\de a,u\quer u\,b]
  =[u\de a,u\quer u\,b]=[\de a,b]$.
\end{bew}

\begin{lemma}\label{bild gleich kern}
Let $(\de a,a),(\de b,b)\in I$ and
$\lambda\de a=\mu\de b$ for some $\lambda,\mu\in\Z[1/n][\xi]$,
$\lambda,\mu\neq 0$. Then $u\in\Z[1/n][\xi +\xi^{-1}]^*$ exists such that
$[\de a,a]=[\de b,ub]$.
\end{lemma}

\begin{bew}
If $\lambda\de a=\mu\de b$, then
$\quer\lambda\quer{\de a}=\quer\mu\quer{\de b}$
and herewith
$$(\lambda\quer\lambda \,a)
  =\lambda\de a\quer\lambda\quer{\de a}
  =\mu\de b\quer\mu\quer{\de b}
  =(\mu\quer\mu \,b).$$
But then a unit $u\in\Z[1/n][\xi +\xi^{-1}]^*$ exists with
$\lambda\quer\lambda \,a=\mu\quer\mu \,ub$. Herewith
$$[\de a,a]
  =[\lambda\de a,\lambda\quer\lambda \,a]
  =[\mu\de b,\mu\quer\mu \,ub]
  =[\de b,ub].$$
\end{bew}

\begin{prop}\label{Struktur der Paare}
Let $\mathcal C_0$ be the ideal class group of $\Z[1/n][\xi]$. Then
the sequence
$$
1
\longrightarrow
\Z[1/n][\xi +\xi^{-1}]^*/N(\Z[1/n][\xi]^*)
\overset{\delta}{\longrightarrow}
\mathcal I
\overset{\eta}{\longrightarrow}
\mathcal C_0
\longrightarrow
1
$$
where $\delta([u])=[\Z[1/n][\xi],u]$, $\eta([\de a,a])=[\de a]$,
is a short exact sequence.
\end{prop}

\begin{bew}
Lemma~\ref{delta injektiv} implies that $\delta$ is injective and $\eta$
is well-defined and surjective. Moreover
$$\eta(\delta([u]))=\eta([\Z[1/n][\xi],u])=[\Z[1/n][\xi]]$$
and Lemma~\ref{bild gleich kern} implies that the kernel of $\eta$
is equal to the image of $\delta$.
\end{bew}

\begin{cor}\label{Anzahl Konj'kl}
The number of conjugacy classes of matrices of order $p$ in
$\SP(p-1,\Z[1/n])$ is equal to
$$|\mathcal C_0|\cdot\bigl[\Z[1/n][\xi +\xi^{-1}]^* : N(\Z[1/n][\xi]^*)\bigr].$$
\end{cor}

\begin{bew}
This corollary is a direct consequence of
Proposition~\ref{Struktur der Paare} because the
number of conjugacy classes of matrices of order $p$ in
$\SP(p-1,\Z[1/n])$ is equal to the cardinality of $\mathcal I$.
\end{bew}

If $\Z[1/n][\xi]$ is a principal ideal domain the cardinality
of $\mathcal C_0$ is $1$ and the number of conjugacy classes of
matrices of order $p$ in $\SP(p-1,\Z[1/n])$ is given only by the index
defined above. In fact we can choose $n\in\Z$ such that $\Z[1/n][\xi]$
is a principal ideal domain. Indeed let
$\de a_1,\dots,\de a_h$ be representatives of the ideal
classes of $\Q(\xi)$. For $j=1,\dots,h$ choose $n_j\in\de a_j$
with $n_j\in\Z[1/n][\xi]$. It is possible to choose the $n_j$ such that
$n_j\in\Z$. Then $n=\prod_{j=1}^h n_j\in\de a_k$ for any $k$ with
$1\leq k\leq h$. For more details see Lang~\cite{langant} and
Neukirch~\cite{neukirch}.

\subsection{The number of conjugacy classes}

Let $N:\Q(\xi)\longrightarrow \Q(\xi +\xi^{-1})$ be the norm mapping
defined above. Let $n\in\Z$ and $\xi$ a primitive $p$th root of unity.
This is the aim of this section is to compute the number of conjugacy
classes of elements of order $p$ in $\SP(p-1,\Z[1/n])$. Therefore
we use Corollary~\ref{Anzahl Konj'kl}.

Kummer proved that
$\Z[1/n][\xi]^* = \Z[1/n][\xi +\xi^{-1}]^*\times \langle -\xi\rangle$ where
$\langle -\xi\rangle$ is the group of roots of unity in $\Q(\xi)$.
This implies that
$$\bigl[\Z[\xi +\xi^{-1}]^* : N(\Z[\xi]^*)\bigr]
  =\bigl[\Z[\xi +\xi^{-1}]^*:(\Z[\xi +\xi^{-1}]^*)^2\bigr].$$
Moreover
$\Z[\xi +\xi^{-1}]^*\cong\Z^{(p-3)/2}\times\Z/2\Z$ because of the Dirichlet
unit theorem. Therefore
$$\bigl[\Z[\xi +\xi^{-1}]^* : N(\Z[\xi]^*)\bigr]=2^{\frac{p-1}{2}}.$$
Since the prime above $p$ in $\Z[\xi]$ is principal, generated by $1-\xi$,
and the prime above $p$ in $\Z[\xi +\xi^{-1}]$ is principal, generated
by $N(1-\xi)=(1-\xi)(1-\xi^{-1})$, we get
$$\bigl[\Z[1/p][\xi +\xi^{-1}]^* : N(\Z[1/p][\xi]^*)\bigr]=2^{\frac{p-1}{2}}.$$

\begin{prop}\label{index bzgl S}
Let $p$ be an odd prime and let $\xi$ be a primitive $p$th root of
unity. Let $S^+$ be a finite set of prime ideals in $\Z[\xi +\xi^{-1}]$,
and let $S$ be the set of the prime ideals in $\Z[\xi]$ that lie over
those in $S^+$. Then
$$\bigl[(\Z[\xi +\xi^{-1}]^{S^+})^*:N((\Z[\xi]^S)^*)\bigr]
   = 2^{\frac{p-1}{2}+\tau}$$
where $\tau$ is the number of inert primes in $S^+$.
\end{prop}

\begin{bew}
Let $S:=\{\de q_1,\dots,\de q_k\}$ be a set of prime ideals in $\Z[\xi]$.
Then the isomorphism given by the generalization of the Dirichlet unit
theorem implies
that for each prime ideal $\de q_j\in S$, $j=1,\dots,k$,
$g_j\in\de q_j$ exists such that each unit $u\in(\Z[\xi]^S)^*$
can be written
$$u=u' g_1^{n_1}\cdots g_k^{n_k}$$
where $u'\in\Z[\xi]^*$, $n_j\in\Z$, $j=1,\dots,k$.
We compute the index we want to know by induction on the number
of primes in $S^+$.
Let $T^+$ be a finite set of prime ideals in $\Z[\xi +\xi^{-1}]$.
Let $T$ be the
set of those prime ideals in $\Z[\xi]$ that lie over the prime
ideals in $T^+$. Define $S^+:=T^+ \cup \{\de q^+\}$ where
$\de q^+ \subset\Z[\xi +\xi^{-1}]$, $\de q^+\not\in T^+$, is a prime ideal.
Let $S$ be the set of the prime ideals in $\Z[\xi]$ that lie over
the prime ideals in $S^+$.
We have the following possibilities.
\begin{enumerate}
\item The prime $\de q^+$ is inert. Then $S=T\cup\{\de q\}$ where
$\de q$ is the prime that lies over $\de q^+$.

\item The prime $\de q^+$ splits in $\Z[\xi]$. Then
$S=T\cup\{\de q,\quer{\de q}\}$
where $\de q$, $\quer{\de q}$ are the primes that lie over $\de q^+$.

\item The prime $\de q^+$ lies over $p$. Then $S=T\cup\{\de p\}$
where $\de p=(1-\xi)$, the prime over $p$.
\end{enumerate}
We have
$$(\Z[\xi +\xi^{-1}]^{S^+})^* \cong (\Z[\xi +\xi^{-1}]^{T^+})^* \times \Z\,.$$
If the prime $\de q^+$ is inert or if it lies over $p$,
cases i) and iii) above, then
\begin{align*}
(\Z[\xi]^S)^* &\cong \Z[\xi]^* \times \Z^{|S|}
               \cong \Z[\xi]^* \times \Z^{|T|} \times \Z \\
              &\cong (\Z[\xi]^T)^* \times \Z
\end{align*}
and if the prime $\de q^+$ splits in $\Z[\xi]$, case ii) above, then
$$(\Z[\xi]^S)^* \cong (\Z[\xi]^T)^* \times \Z^2.$$
We give a formula for the index
$$\bigl[(\Z[\xi +\xi^{-1}]^{S^+})^*:N((\Z[\xi]^S)^*)\bigr]$$
in relation to the index
$$\bigl[(\Z[\xi +\xi^{-1}]^{T^+})^*:N((\Z[\xi]^T)^*)\bigr].$$
If the prime $\de q^+$ is inert, then
$$\bigl[(\Z[\xi +\xi^{-1}]^{S^+})^*:N((\Z[\xi]^S)^*)\bigr] =
   2\bigl[(\Z[\xi +\xi^{-1}]^{T^+})^*:N((\Z[\xi]^T)^*)\bigr].$$
If the prime $\de q^+$ splits in $\Z[\xi]$ or if it lies over $p$, then
$$\bigl[(\Z[\xi +\xi^{-1}]^{S^+})^*:N((\Z[\xi]^S)^*)\bigr] =
   \bigl[(\Z[\xi +\xi^{-1}]^{T^+})^*:N((\Z[\xi]^T)^*)\bigr].$$
This shows that if we add an inert prime to the set $S$ the index
is multiplied by $2$, and if we add primes that split or the prime
over $p$, then the index does not change.
\end{bew}

\begin{satz}\label{Index ist}
Let $n\in\Z$. Then
$$\left[\Z[1/n][\xi +\xi^{-1}]^*:N(\Z[1/n][\xi]^*)\right]
=2^{\frac{p-1}{2}+\tau}$$
where $\tau$ is the number of inert primes in
$\Z[\xi +\xi^{-1}]$ that lie over primes in $\Z$ that divide $n$.
\end{satz}

\begin{bew}
Let $n\in\Z$ and let $S^+$, resp.\ $S$, be the prime ideals in
$\Z[\xi +\xi^{-1}]$, resp.\ $\Z[\xi]$,
over the primes in $\Z$ that divide $n$.
Then the assumption
follows directly from Proposition~\ref{index bzgl S}.
\end{bew}

\begin{satz}\label{Konj'kl von Elem. in Sp(p-1,Z[1/n])}
The number of conjugacy classes of matrices of order $p$ in
$\SP(p-1,\Z[1/n])$, $0\neq n\in\Z$, is
$$|\mathcal C_0| 2^{\frac{p-1}{2}+\tau}$$
where $\mathcal C_0$ is the ideal class group of $\Z[1/n][\xi]$ and
$\tau$ is the number of inert primes in $\Z[\xi +\xi^{-1}]$ that lie over
primes in $\Z$ that divide $n$.
\end{satz}

\begin{bew}
This follows directly from Corollary~\ref{Anzahl Konj'kl} and
Theorem~\ref{Index ist}.
\end{bew}


\section{Subgroups of order $p$}


\subsection{The quotient of the normalizer by the centralizer of
subgroups of order $p$}

The aim is to study the centralizers and normalizers of conjugacy
classes of subgroups of order $p$ in $\SP(p-1,\Z[1/n])$. We use
the bijection between the set $\mathcal I$ of equivalence classes
$[\de a,a]$ and the conjugacy classes of matrices of order $p$.
Each conjugacy class of matrices generates a conjugacy class of
subgroups of order $p$ in $\SP(p-1,\Z[1/n])$. We determine the
equivalence classes $[\de a,a]$ that correspond to the conjugacy
classes of the elements of a subgroup.

Let $Y\in\SP(p-1,\Z[1/n])$ be of odd prime order $p$.
We have seen that the conjugacy class of $Y$ corresponds to an
equivalence class $[\de a,a]$. Let
$$\alpha
  = (\alpha_1,\dots,\alpha_{p-1})\tp
  \in \bigl(\Z[1/n][\xi]\bigr)^{p-1}$$
be an eigenvector of $Y$ to the eigenvalue $\xi=e^{i2\pi/p}$.
It is obvious that $Y^l=\xi^l\alpha$. Let
$\gamma_k\in\Gal(\Q[\xi]/\Q)$ be such that $\gamma_k(\xi)=\xi^k$.
Then $\gamma_k(\xi^l)=\xi^{kl}$. If $kl\equiv 1\mod p$, then
$\gamma_k(\xi^l)=\xi$ and moreover
$$
Y^l\gamma_k(\alpha)
= \gamma_k(Y^l\alpha) = \gamma_k(\xi^l\alpha)
   = \gamma_k(\xi^l)\gamma_k(\alpha) = \xi^{kl}\gamma_k(\alpha)
= \xi\gamma_k(\alpha).
$$
So $\gamma_k(\alpha)$ is the eigenvector of $Y^l$ to the eigenvalue $\xi$.
Let $\de b$ be the ideal given by the $\Z[1/n]$-basis
$\gamma_k(\alpha_1),\dots,\gamma_k(\alpha_{p-1})$. Moreover
let
$$b = D^{-1}(\gamma_k(\alpha))\tp J\gamma_k(\quer\alpha)
    = D^{-1}\gamma_k(\alpha\tp J\quer\alpha).$$
So the conjugacy class of $Y^l$ corresponds to the
equivalence class $[\de b,b]$ with
\begin{align*}
\de b &= \gamma_k(\de a)\\
b &= D^{-1}\gamma_k(Da) = D^{-1}\gamma_k(D)\gamma_k(a).
\end{align*}
Let $S$ be a multiplicative set such that $S^{-1}\Z =\Z[1/n]$.
Then $S^{-1}\Z[\xi] =\Z[1/n][\xi]$ and the different in
$\Z[\xi]$ and in $\Z[1/n][\xi]$ are both principal
ideals generated by $D=p\,\xi^{(p+1)/2}/(\xi -1)$.
If $p\mid n$, then $D$ is a unit in $\Z[1/n][\xi]$ since $(\xi -1)$
is a prime that divides $p$. If $u,v\in\Z[1/n][\xi]^*$
are units with $u=\quer u$, $v=-\quer v$, then $Du=-\quer{Du}$
and $Dv=\quer{Dv}$. This shows that the multiplication with $D$
defines an isomorphism on $\Z[1/n][\xi]^*$ that yields a
bijection between the real and the purely imaginary units.

\begin{satz}\label{N/C in Sp(p-1,Z[1/n])}
Let $n\in\Z$ be such that $\Z[1/n][\xi]$ and $\Z[1/n][\xi +\xi^{-1}]$
are principal ideal domains and moreover $p\mid n$. Let $N(P)$ denote
the normalizer and $C(P)$ the centralizer of a subgroup $P$ of order
$p$ in $\SP(p-1,\Z[1/n])$. Then
$$N(P)/C(P)\cong\Z/j\Z$$
where $j\mid p-1$, $j$ odd. For each $j$ with $j\mid p-1$, $j$ odd,
exists a subgroup of order $p$ in $\SP(p-1,\Z[1/n])$
with $N(P)/C(P)\cong\Z/j\Z$.
\end{satz}

\begin{bew}
Let $n$ be such that $\Z[1/n][\xi]$ is a principal ideal domain.
Then the ideal $\de a$ in the pair $[\de a,a]$ is a principal ideal.
If $\de a=(x)$, then $(x)\quer{(x)}=(x\quer x)=(a)$,
i.e.\ a unit $u$ exists such that $a=ux\quer x$. Then
$$[\de a,a]=[(x),a]=[\Z[1/n][\xi],u].$$
The conjugacy class of $Y\in\SP(p-1,\Z[1/n])$ corresponds
to $[\Z[1/n][\xi],u]$. For $1<l<p-1$ we have seen that the
conjugacy class of $Y^l$ corresponds to
$[\Z[1/n][\xi],D^{-1}\gamma_k(Du)]$ where
$\gamma_k\in\Gal(\Q[\xi]/\Q)$ is defined such that
$\gamma_k(\xi^l)=\xi$. The matrices $Y$ and $Y^l$ are
conjugate if and only if
$$[\Z[1/n][\xi],u]=[\Z[1/n][\xi],D^{-1}\gamma_k(Du)].$$
Lemma~\ref{delta injektiv} shows that this equation is
satisfied if and only if $\omega\in\Z[1/n][\xi]^*$ exists
such that
\begin{equation}\label{stern}
D^{-1}\gamma_k(Du)=u\omega\quer\omega .
\end{equation}
We know that $u\in\Z[1/n][\xi +\xi^{-1}]^*$ and this implies that $Du$
is purely imaginary. First we check if $u\in\Z[1/n][\xi +\xi^{-1}]^*$
exists such that a special case of (\ref{stern}) holds, namely
the case with $\omega=1$, i.e.\ we try to find $\gamma_k$ and
$u$ such that $\gamma_k(Du)=Du$. The automorphism
$\gamma_{p-1}(=\gamma_{-1})$ has order $2$, i.e.\ $\gamma_{p-1}$
yields the complex conjugation.
Since $u$ is real and therefore $Du$ purely imaginary,
we get $\gamma_{p-1}(Du)=-Du$. This proves that neither
$\gamma_k(Du)=Du$ nor (\ref{stern}) can be satisfied
if $k=p-1$ (the image of $\omega\quer\omega$ under any
embedding of $\Z[1/n][\xi]$ in $\C$ is a positive real
number).  Any automorphism $\gamma_k\in\Gal(\Q[\xi]/\Q)$
generates a subgroup
$\langle\gamma_k\rangle\subseteq\Gal(\Q[\xi]/\Q)$ and the
order of this subgroup divides $p-1$, the order of
$\Gal(\Q[\xi]/\Q)$. Let $j=|\langle\gamma_k\rangle|$
denote the order of $\gamma_k$. If $j$ is even
the order of $\gamma_k^{j/2}$ is $2$ and on the other hand
$\gamma_k^r(Du)=Du$ for any $1<r<j$. This yields a
contradiction and therefore $\gamma_k(Du)=Du$ cannot
be satisfied if the order of $\gamma_k$ is even.
This implies that if $\gamma_k$ and $u$ exist with
$\gamma_k(Du)=Du$, then the order of $\gamma_k$ is odd.

The main theorem of Galois theory says that a subfield
$\Q\subseteq K\subseteq\Q(\xi)$ corresponds to the subgroup
$\langle\gamma_k\rangle\subseteq\Gal(\Q[\xi]/\Q)$ and that
$$K=\{x\in\Q(\xi)\mid
     \forall\gamma_{k^r}\in\langle\gamma_k\rangle,\
     \gamma_{k^r}(x)=x \}.$$
Let $n\in\Z$ with $p\mid n$, We have seen that in this case
$D=p\xi^{(p+1)/2}/(\xi -1)$ is a unit in $\Z[1/n][\xi]$.
We also know that $D=-\quer D$. Let
$\gamma_k\in\Gal(\Q[\xi]/\Q)$ be of odd order $j$. Since
complex conjugation commutes with the Galois automorphisms,
we get for any $r$, $1\leq r\leq j$,
$\gamma_{k^r}(D)=-\quer{\gamma_{k^r}(D)}$. Since $j$ is odd,
$$\prod_{r=1}^j\gamma_{k^r}(D)
  =(-1)^j\prod_{r=1}^j\gamma_{k^r}(\quer D)
  = -\prod_{r=1}^j\gamma_{k^r}(\quer D).$$
Moreover this product is invariant under $\gamma_k$ since
$$\gamma_k\Bigl(\prod_{r=1}^j\gamma_{k^r}(D)\Bigr)
  =\prod_{r=1}^j\gamma_k\bigl(\gamma_{k^r}(D)\bigr)
  =\prod_{r=1}^j\gamma_{k^r}(D).$$
Now consider the composition
$\gamma_k\circ\gamma_{p-1}=\gamma_{-k}$ where the order
of $\gamma_k$ is odd. The order of $\gamma_{-k}$ is
even and $\langle\gamma_k\rangle$ is a subgroup of
$\langle\gamma_{-k}\rangle$. Let $L$ denote the
subfield $\Q\subseteq L\subseteq K\subseteq\Q(\xi)$
corresponding to $\langle\gamma_{-k}\rangle$.
Sinnott constructs in \cite{sinnott2} cyclotomic units in
any subfield $L$ of $\Q(\xi_m)$ where $\xi_m$ is a $m$th
root of unity. This means that units exist in $L$, that
are contained in no subfield of $L$. Let $v\in L$ be such a
unit ($v\in\Z[\xi]$). Then $\gamma_{p-1}(v)=v$,
$\gamma_k(v)=v$ since $\langle\gamma_{-k}\rangle$ fixes
the elements of $L$.
Let $w:=\prod_{r=1}^{j-1}\gamma_{k^r}(D)\in\Z[1/n][\xi]$.
Then
$$
w = \prod_{r=1}^{j-1}\gamma_{k^r}(D)
   = (-1)^{j-1}\prod_{r=1}^{j-1}\quer{\gamma_{k^r}(D)}
   = \quer w
$$
since $j$ is odd.
Moreover $Dw=\prod_{r=1}^j\gamma_{k^r}(D)$ and therefore
$\gamma_k(Dw)=Dw$. We have $Dw=-\quer{Dw}$ since $w=\quer w$.
Now $wv=\quer{wv}$ is a unit. Let $u=wv$, then the construction implies
that
$$\gamma_k(Du)=\gamma_k(Dwv)=Dwv=Du.$$
So for any $\gamma_k\in\Gal(\Q(\xi)/\Q)$ of odd order, we found
$u\in\Z[1/n][\xi +\xi^{-1}]^*$ with $\gamma_k(Du)=Du$ and such that
$$[\Z[1/n][\xi],u]=[\Z[1/n][\xi],D^{-1}\gamma_k(Du)].$$
If $Y\in\SP(p-1,\Z[1/n])$ is in the corresponding equivalence
class then this is also true for $Y^l$ with $l$ such that
$\gamma_k(\xi^l)=\xi$. If $Y$ is conjugate to $Y^l$ with
$\gamma_k(\xi^l)=\xi$, then $Y$ is also conjugate to
$Y^{l^r}$ where $1\leq r\leq j$ and $j$ is the order of
$\gamma_k$. Indeed $\gamma_{k^r}(\xi^{l^r})=\xi$ for
$1\leq r\leq j$ and therefore $l^j\equiv 1\mod p$ (since
$\gamma_{k^j}=id$) and $Y^{l^j}=Y$ because the order of $Y$
is $p$. The $l^r$ form a cyclic subgroup of $\Z/p\Z$.

Let $k,l\in\Z$ be as above, i.e. $\gamma_k(\xi^l)=\xi$.
Let $\gamma_l\in\Gal(\Q(\xi)/\Q)$ with $\gamma_l(\xi)=\xi^l$.
Then $\gamma_l=\gamma_k^{-1}$ and if $j$ is the order of
$\gamma_k$, then $j$ is also the order of $\gamma_l$.
Therefore $l^j\equiv 1\mod p$. This means that $Y^{l^j}=Y$
and the $Y^{l^r}$, $1\leq r\leq j$ are conjugate to $Y$.
We know that $j$ is odd and $j\mid p-1$.

If $j$ elements are conjugate in the subgroup generated by the
matrix $Y\in\SP(p-1,\Z[1/n])$, and if $j$ is maximal with this
property, then we have for this subgroup $N(P)/C(P)\cong\Z/j\Z$
since $\langle\gamma_k\rangle\cong\Z/j\Z$.
Since we showed that for any odd divisor $j\mid p-1$ a matrix
$Y\in\SP(p-1,\Z[1/n])$ exists for which $j$ powers are
conjugate, we showed that for any $j\mid p-1$, $j$ odd, a
subgroup of order $p$ exists in $\SP(p-1,\Z[1/n])$, for
which $N(P)/C(P)\cong\Z/j\Z$.
\end{bew}

\subsection{The centralizer of subgroups of order $p$}

\begin{satz}\label{C in Sp(p-1,Z[1/n])}
Let $n\in\Z$ be such that $\Z[1/n][\xi]$ and $\Z[1/n][\xi +\xi^{-1}]$ are
principal ideal domains. Then the centralizer $C(P)$ of a
subgroup $P$ of order $p$ in $\SP(p-1,\Z[1/n])$ is
$$C(P)\cong \Z/2p\Z\times\Z^{\sigma^+}$$
where $\sigma^+ =\sigma$ if $p\nmid n$, $\sigma^+=\sigma +1$ if
$p\mid n$ and $\sigma$ is the number of primes in $\Z[\xi +\xi^{-1}]$ that
split in $\Z[\xi]$ and lie over primes in $\Z$ that divide $n$.
\end{satz}

\begin{bew}
Let $Y\in\SP(p-1,\Z[1/n])$ be of order $p$ and let
$[\de a,a]$ be the equivalence class corresponding
to the conjugacy class of $Y$.
Let $P$ be the subgroup generated by $Y$. Let
$Z\in\SP(p-1,\Z[1/n])$ be an element of the centralizer
of $Y$, i.e. $Z^{-1}YZ=Y$ or $YZ=ZY$. Then $Z$ is an element
of the centralizer of $P$. If $\alpha$ is an
eigenvector of $Y$ to the eigenvalue $\xi$, then so is $Z\alpha$:
$$\xi Z\alpha =Z\xi\alpha =ZY\alpha =YZ\alpha.$$
But this means that $Z\alpha=w\alpha$
for some $w\in\Z[1/n][\xi]$ and $w$ is a unit since $Z$
is invertible. Therefore
\begin{align*}
(Z\alpha)\tp J\quer{Z\alpha}^{(i)}
&= \alpha\tp Z\tp JZ\quer\alpha^{(i)}
 = w\alpha\tp J\quer w^{(i)}\quer\alpha^{(i)}
\\
&= w\quer w^{(i)}\alpha\tp J\quer\alpha^{(i)}
 = \delta_{1i}aw\quer w^{(i)}D
\end{align*}
and, since $\delta_{1i}=0$ for $i\neq 1$, we get
$$(Z\alpha)\tp J\quer{Z\alpha}=aw\quer wD.$$
But $Z\in\SP(p-1,\Z[1/n])$ and therefore
$$
(Z\alpha)\tp J\quer{Z\alpha}
 = \alpha\tp Z\tp JZ\quer\alpha
 = \alpha\tp J\quer\alpha
 = aD.
$$
This implies that $w\quer w=1$.
In order to determine the centralizer $C(P)$ of a subgroup
$P\subseteq \SP(p-1,\Z[1/n])$ of order $p$, we have to
find the units $w\in\Z[1/n][\xi]^*$ that satisfy $w\quer w =1$.
This corresponds to the kernel of the norm mapping
$$\begin{array}{rccl}
N: & \Z[1/n][\xi]^* & \longrightarrow & \Z[1/n][\xi +\xi^{-1}]^* \\
 & x & \longmapsto & x\quer x.
\end{array}$$
Brown \cite{browna} and Sjerve and Yang \cite{sjerve} showed
that the kernel of the norm mapping
$$\begin{array}{rccl}
N': & \Z[\xi]^* & \longrightarrow & \Z[\xi +\xi^{-1}]^* \\
 & x & \longmapsto & x\quer x
\end{array}$$
is the set of roots of unity
$$\ker(N')=\{\pm\xi^r\mid \xi^p=1, 1\leq r\leq p\}.$$
It is obvious that $\ker(N')\subseteq\ker(N)$.
The prime ideals that lie over the primes in $\Z$ and divide $n$ yield
units in $\Z[1/n][\xi]^*\setminus\Z[\xi]^*$. Let
$\de q^+\subseteq\Z[\xi +\xi^{-1}]$ be a prime over a prime $q\mid n$ and let
$\de q\subseteq\Z[\xi]$ be a prime over $\de q^+$. If $\de q^+$ is
inert, then $\de q=\quer{\de q}$ and if $\de q^+$ splits, then
$\de q^+\Z[\xi] =\de q\quer{\de q}$. A generalization for $S$-units of the
Dirichlet unit theorem says that for each  prime $\de q_j$,
$j=1,\dots,k$, over $n$ a $g_j\in\de q_j$ exists such that any unit
$u\in(\Z[1/n][\xi])^*$ can be written as
$$u=u' g_1^{n_1}\cdots g_k^{n_k}$$
where $u'\in\Z[\xi]^*$, $n_j\in\Z$, $j=1,\dots,k$.
So the group of units $\Z[1/n][\xi +\xi^{-1}]^*$
is generated by $\Z[\xi +\xi^{-1}]^*$, the inert primes over $n$, the
primes over $n$ that split and, if $p\mid n$, the prime over $p$.
The inert primes yield nontrivial elements in
$\Z[1/n][\xi +\xi^{-1}]^*/N(\Z[1/n][\xi]^*)$ since for those holds
$w\quer w=w^2\neq 1$ for $w\neq\pm 1$.
The centralizer $C(P)$ is a finitely generated group whose
torsion subgroup is isomorphic to the group of roots of unity
in $\Q(\xi)$ and whose rank is equal to $\sigma$ if $p\nmid n$ and
to $\sigma + 1$ if $p\mid n$ where
$$\sigma^+ = \operatorname{rank}\bigl(\Z[1/n][\xi]^*\bigr)
          -\operatorname{rank}\bigl(\Z[1/n][\xi +\xi^{-1}]^*\bigr).$$
This difference is equal to the number of primes in $\Z[\xi +\xi^{-1}]$ that
split or ramify in $\Z[\xi]$ and lie over primes in $\Z$ that divide $n$.
This follows directly from a generalization of the Dirichlet unit
theorem and proves our theorem.
\end{bew}

\subsection{The action of the normalizer on the centralizer of
 subgroups of order $p$}

\begin{satz}\label{Wirkung von N auf C}
Let $N(P)$ be the normalizer and $C(P)$ the centralizer of a subgroup
$P$ of order $p$ in $\SP(p-1,\Z[1/n])$. Let $p$ be an odd prime,
$\xi$ a primitive $p$th root of unity, $n\in\Z$ such that
$\Z[1/n][\xi]$ and $\Z[1/n][\xi +\xi^{-1}]$ are principal ideal domains
and moreover $p\mid n$.
Then the action of $N(P)/C(P)$ on $C(P)$ is given by the action of the
Galois group $\Gal(\Q(\xi)/\Q)$ on the group of units
$\Z[1/n][\xi]^*$. Moreover $N(P)/C(P)$ acts faithfully on $C(P)$.
\end{satz}

\begin{bew}
We have seen in the proof of
Theorem~\ref{C in Sp(p-1,Z[1/n])} that the centralizer of a subgroup
of order $p$ in $\SP(p-1,\Z[1/n])$ is given by the kernel of the
norm mapping $\Z[1/n][\xi]^* \longrightarrow \Z[1/n][\xi +\xi^{-1}]^*$,
$x \longmapsto x\quer x$. Herewith the centralizer is isomorphic
to a subgroup of the group of units $\Z[1/n][\xi]^*$. In the proof of
Theorem~\ref{N/C in Sp(p-1,Z[1/n])} we identify the quotient
$N(P)/C(P)$ with a subgroup of the Galois group $\Gal(\Q(\xi)/\Q)$.
Herewith the action of the quotient $N(P)/C(P)$ on the centralizer
$C(P)$ is given by the action of the subgroup of $\Gal(\Q(\xi)/\Q)$
corresponding to $N(P)/C(P)$ on the kernel of the norm mapping
$\Z[1/n][\xi]^* \longrightarrow \Z[1/n][\xi +\xi^{-1}]^*$. Since it is
nontrivial, the action of $N(P)/C(P)$ on $C(P)$ is faithful.
\end{bew}

\vspace*{1cm}

\begin{minipage}[t]{12cm}
Cornelia M. Busch

Katholische Universit\"at Eichst\"att-Ingolstadt

MGF

D-85071 Eichst\"att

Germany
\end{minipage}


\begin{thebibliography}{99}
\bibitem{browna}
        K.\ S.\ Brown,
        \emph{Euler {C}haracteristics of {D}iscrete {G}roups and ${G}$-{S}paces},
        Invent. Math. 27 (1974), 229-264.

\bibitem{brownb}
        K.\ S.\ Brown,
        \emph{Cohomology of Groups}, GTM 87, Springer 1982.

\bibitem{BuschSC}
        C.\ Busch,
        \emph{Symplectic characteristic classes},
        L'Enseignement Math\'e\-matique 47 (2001), 115-130.

\bibitem{BuschHFSp}
        C.\ Busch,
        \emph{The Farrell cohomology of \hbox{$\mathrm{Sp}(p-1,\mathbb{Z})$}},
        Documenta Mathematica 7 (2002), 239-254.

\bibitem{buschperFar}
        C.~M. Busch,
        \emph{On $p$-periodicity in the Farrell cohomology of
        \hbox{$\mathrm{Sp}(p-1,\mathbb{Z}[1/n])$}}, Preprint (2005).

\bibitem{langant}
        S. Lang,
        \emph{Algebraic number theory},
        Addison Wesley 1970.

\bibitem{nadim}
        N. Naffah,
        \emph{On the {I}ntegral {F}arrell {C}ohomology {R}ing of ${PSL}_2(\mathbb{Z}[1/n])$},
        Diss. {E}{T}{H} {N}o. 11675, ETH Z\"urich, 1996.

\bibitem{neukirch}
        J.\ Neukirch,
        \emph{Algebraic number theory},
        Grundlehren der mathematischen Wissenschaften 322, Springer 1999.

\bibitem{serrelf}
        J.-P. Serre,
        \emph{Local Fields}, GTM 67, Springer 1979.

\bibitem{sinnott2}
        W. Sinnott,
        \emph{On the {S}tickelberger ideal and the circular units of an abelian field},
        Invent. Math. 62 (1980), 181-234.

\bibitem{sjerve}
        D. Sjerve and Q. Yang,
        \emph{Conjugacy {C}lasses of $p$-{T}orsion in $\mathrm{Sp}_{p-1}(\mathbb{Z})$},
        J. of Algebra 195 (1997), 580-603.

\bibitem{wash}
        L. C. Washington
        \emph{Introduction to cyclotomic fields}, GTM 83, Springer 1997.
\end{thebibliography}
\end{document}